\documentclass[oneside,12pt]{amsart}

\usepackage{amscd,amsmath,latexsym}
\usepackage[mathcal]{euscript}

\newtheorem{thm}{Theorem}[section]
\newtheorem{cor}[thm]{Corollary}
\newtheorem{lem}[thm]{Lemma}
\newtheorem{prop}[thm]{Proposition}
\theoremstyle{definition}
\newtheorem{defin}[thm]{Definition}
\newtheorem{exm}[thm]{Example}

\theoremstyle{remark}
\newtheorem*{rem}{Remark}

\newcommand{\Z}{{\mathbb Z}}

\begin{document}

\title{Commuting elements and spaces of homomorphisms}

\author[Alejandro ~Adem]{Alejandro ~Adem$^{*}$}
\address{Department of Mathematics
University of British Columbia, Vancouver, B.~C., Canada}
\email{adem@math.ubc.ca}
\thanks{$^{*}$Partially supported by the NSF and NSERC}

\author[Frederick ~R.~Cohen]{Frederick ~R.~Cohen$^{**}$}
\address{Department of Mathematics,
University of Rochester, Rochester, NY 14627}
\email{cohf@math.rochester.edu}
\thanks{$^{**}$Partially supported by the NSF}
\date{\today}

\begin{abstract}

This article records basic topological, as well as homological
properties of the space of homomorphisms $Hom(\pi,G)$ where $\pi$ is
a finitely generated discrete group, and $G$ is a Lie group,
possibly non-compact. If $\pi$ is a free abelian group of rank equal
to $n$, then $Hom(\pi, G)$ is the space of ordered $n$--tuples of
commuting elements in $G$. If $G=SU(2)$, a complete calculation of
the cohomology of these spaces is given for $n=2, 3$. An explicit
stable splitting of these spaces is also obtained, as a special case
of a more general splitting.
\end{abstract}

\maketitle

\section{Introduction}
\label{sec:Introduction}

Let $\pi$ denote a finitely generated discrete group and $G$ a Lie
group. The set of homomorphisms $Hom(\pi, G)$ from $\pi$ to $G$ can
be topologized in a natural way as a subspace of a finite product of
copies of $G$, where the relations on $\pi$ carve out what can be a
variety with interesting singularities. The purpose of this paper is
to describe the geometry and cohomology of these spaces for certain
classes of groups.

These spaces are interesting and can be quite complicated. The
approach in this article is to start by studying their basic
structure. Additional methods for analyzing $Hom(\pi, G)$ are
developed by restricting to suitable classes of discrete groups
$\pi$, and varying the target group $G$.

First consider the question of path--connectivity for these spaces.
In the case of $G=O(n)$ or $G=SO(n)$, the first two Stiefel--Whitney
classes play an important role:

\begin{thm}
For a finitely generated discrete group $\pi$, the following two
statements hold:
\begin{itemize}

\item
There is a decomposition into non--empty, disjoint closed subsets
(called $w$--sectors):
$$Hom(\pi, O(n))\cong \bigsqcup_{w\in H^1(\pi,\mathbb F_2)}
Hom(\pi, O(n))_w$$ where a $w$--sector corresponds to those
homomorphisms with first Stiefel--Whitney class equal to $w$. In
particular if $H^1(\pi,\mathbb F_2)\ne 0$, then $Hom(\pi, O(n))$ is
not path--connected. The sector for $w=0$ is precisely $Hom(\pi,
SO(n))$.

\item
For $n$ sufficiently large, $Hom(\pi, O(n))_w$ has at least as many
components as the cardinality of the image of
$H^2(\pi/[\pi,\pi],\mathbb F_2)\to H^2(\pi, \mathbb F_2)$.
\end{itemize}
\end{thm}

This applies very well to fundamental groups of certain complements
of complex hyperplane arrangements, such as the pure braid groups:

\begin{cor}
If $\pi$ is the fundamental group of the complement of an
arrangement of complex hyperplanes which has contractible universal
cover, and if $n$ is sufficiently large, then
$$\# \pi_0(Hom(\pi, O(n))\ge
|H^1(\pi,\mathbb F_2)||H^2(\pi,\mathbb F_2)|$$
\end{cor}

Kapovich and Millson have investigated the geometry of $Hom(\pi, G)$
near the identity for $\pi$ equal to the fundamental group of
complements of certain arrangements \cite{KapMill}, Theorem 1.13.
The approach in this article is to give cruder, more global
information such as the number of path-components.

The number of path-components may change when the target is
restricted from a non--compact Lie group $G$ to a maximal compact
subgroup $K$. Examples are given in \S 3 where the natural map
$Hom(\pi, K) \to Hom(\pi, G)$ fails to be surjective on
path-components. Classically, these arise in the case where $\pi$ is
the fundamental group of a closed orientable surface of genus at
least two with $G = SL(2,\mathbb R)$ and $K = SO(2)$. A further
example given here arises from $\pi = \Gamma_2$ the mapping class
group of a surface of genus two, where $U(2)\subset Sp_4(\mathbb R)$
is a maximal compact subgroup. In this case, $Hom(\Gamma_2,U(2))\to
Hom(\Gamma_2, Sp_4(\mathbb R))$ fails to be surjective on path
components (see \ref{thm:maximal.compact.in.G}).

Under certain conditions, the space of homomorphisms \textsl{must}
be path-connected. This is the case when $\pi=\mathbb Z^n$ and $G$
is a Lie group with path--connected maximal abelian subgroups. For
certain classes of groups this can then be used to analyze the
natural map $Hom(\pi, G)\to [B\pi, BG]$, where $[B\pi, BG]$ denotes
the based homotopy classes of maps from $B\pi$ to $BG$. The
following theorem is proved in section \ref{sec:path components}:

\begin{thm}\label{thm:homologically.toroidal.reps}
Let $\pi$ denote a homologically toroidal group\footnote{A group
$\pi$ is said to be homologically toroidal if there exists a finite
free product $H$ of free abelian groups of finite rank and a
homomorphism $H\to\pi$ which induces a split epimorphism in integral
homology.}; then, if $G$ is a Lie group with path--connected maximal
abelian subgroups, the map
$$B_0:\pi_0(Hom(\pi, G))\to [B\pi, BG]$$
is trivial. In particular, if $\mathcal{M}$ is the complement of an
arrangement of complex hyperplanes which is aspherical, then any
complex unitary representation of $\pi = \pi_1(\mathcal{M})$ induces
a trivial bundle over $B\pi$.
\end{thm}

The second part of this paper is specialized to the cohomology of
the space of ordered commuting $n$--tuples of elements in a compact
Lie group $G$. The main results here mostly apply to the case when
$G=SU(2)$. The approach is to consider the complement, namely the
space of ordered \textsl{non--commuting} $n$--tuples; which seems to
be much more tractable. Duality is then used to compute the
cohomology of the commuting $n$--tuples. For example, the complement
of $Hom(\mathbb Z^2, SU(2))$ in $SU(2)\times SU(2)$ is homotopy
equivalent to $SO(3)$, and the cohomology of the commuting pairs can
be easily computed from this. The case of commuting triples in
$SU(2)^3$ requires more intricate calculations, which yield the
following results:

\begin{thm}\label{thm:cohomology.commuting.triples.SU(2)}

\[ H^i(Hom (\Z\oplus\Z, SU(2)), \Z)
\cong\left\{\begin{array}{r@{\quad\hbox{if}\quad}l}
\Z & i=0\\
0  & i=1\\
\Z & i=2\\
\Z\oplus\Z & i=3\\
\Z/2\Z & i=4\\
0 & i\ge 5\end{array}\right.\]

\[ H^i(Hom (\Z\oplus\Z\oplus\Z, SU(2)), \Z)
\cong\left\{\begin{array}{r@{\quad\hbox{if}\quad}l}
\Z & i=0\\
0  & i=1\\
\Z\oplus\Z\oplus\Z & i=2\\
\Z\oplus\Z \oplus\Z & i=3\\
\Z/2\Z\oplus\Z/2\Z\oplus\Z/2\Z & i=4\\
\Z & i=5\\
0 & i=6\\
\Z/2\Z & i= 7
\\0  & i\ge 8\end{array}\right.\]

\end{thm}

In the last part of this paper, a stable splitting is given for the
space of commuting $n$--tuples. Natural subspaces of $Hom(\mathbb
Z^{n},G)$ arise from the so-called fat wedge filtration of the
product $G^n$ where the base-point of $G$ is $1_G$. Thus if $F_jG^n$
is the subspace of $G^n$ with at least $j$ coordinates equal to
$1_G$, define subspaces of $Hom(\mathbb Z^{n},G)$ by the formula
$S_{n}(j,G) = Hom(\mathbb Z^{n},G)\cap F_jG^n$, and write
$S_n(G)=S_n(1,G)$.

A Lie  group $G$ is said have
\textsl{cofibrantly commuting elements} if the natural inclusions
$I_j: S_{n}(j,G)\to S_{n}(j-1,G)$ are cofibrations for all $n$ and
$j$ for which both spaces are non--empty.  Many Lie groups satisfy
this last property.

\begin{thm} \label{thm:closed.subgroups.of.GL.}
If $G$ is a closed subgroup of $GL(n, \mathbb C)$, then $G$ has
cofibrantly commuting elements.
\end{thm}

In what follows, if $X$ and $Y$ are pointed topological spaces, then
$X\bigvee Y$ denotes their one point union in the product $X\times
Y$ and $\Sigma X$ denotes the suspension of $X$.

\begin{thm} \label{thm:stable.decompositions.for.general.G} If $G$ is a Lie
group which has cofibrantly commuting elements, then there are
homotopy equivalences
$$\Sigma ( Hom(\mathbb Z^n,G))\simeq
\bigvee_{1 \leq k \leq n}\Sigma ( \bigvee^{\binom n k} Hom(\mathbb
Z^k,G)/ S_k(G)).$$
\end{thm}

The last two theorems imply the next corollary.

\begin{cor} \label{cor:stable.homology.decompositions.for.general.G}
If $G$ is a closed subgroup of $GL(n, \mathbb C)$, then there are
isomorphisms of graded abelian groups in cohomology (which may not
preserve products) $$H^*(Hom(\mathbb Z^n,G),\mathbb Z) \to
H^*(\bigvee_{1 \leq k \leq n}\bigvee^{\binom n k} Hom(\mathbb
Z^k,G)/ S_k(G), \mathbb Z).$$
\end{cor}

Applying this to the case of ordered commuting pairs in $SU(2)$,
there is a homotopy equivalence which holds after a single
suspension:

$$Hom( \mathbb Z^2,G)\simeq SU(2)\bigvee SU(2)\bigvee
(\mathbb S^6 - SO(3))$$ where $(\mathbb S^6-SO(3))$ is the
Spanier--Whitehead dual of $SO(3)$ in $\mathbb S^6$. These methods
also yield the result that after one suspension,
$Hom(\mathbb Z^3,SU(2))$ is homotopy equivalent to
$$\bigvee^3 SU(2)\bigvee [\bigvee^3 (\mathbb S^6 - SO(3))] \bigvee
[SU(2) \wedge (\mathbb S^6 -SO(3)].$$ A general decomposition for
$Hom(\mathbb Z^n,G)$ is described, but the summands have not yet
been identified as concretely in terms of Spanier-Whitehead duals.

The authors would like to acknowledge the hospitality and support
of the Max-Planck-Institut f\"ur Mathematik in Bonn, 
the Pacific Institute for the Mathematical Sciences (Vancouver)
and the Institute for Advanced Study (Princeton). They
would like to thank P. Selick and L. Jeffrey for pointing out a
mistake in the homology calculation in an earlier version of this
paper; Enrique Torres for his careful reading of a preliminary
version of this manuscript; and Bill Browder, Karl H. Hofmann and
D. Z. Djokovic for conversations concerning Lie groups.

\section{General Properties of $Hom(\pi, G)$ }
\label{sec:path components}

Basic general properties of spaces of homomorphisms $Hom(\pi, G)$
are given in this section. A key reference here is the paper by
W.Goldman \cite{Goldman}, although the approach will be much less
geometric.

\begin{defin}
Let $\pi$ denote a finitely generated group, and $G$ a finite
dimensional Lie group; $Hom(\pi, G)$ will denote the set of all
group homomorphisms from $\pi$ to $G$.
\end{defin}

Let $F_n$ denote the free group with $n$ generators $\{x_1,x_2,
\cdots, x_n\}$. There is an evaluation map $E:Hom(F_n, G) \to G^n$
defined by $E(f) = (f(x_1),f(x_2), \cdots, f(x_n)).$ There is a
second evaluation map $e: G^n \to Hom(F_n, G)$ defined by
$e(g_1,g_2, \cdots, g_n) = f$ where $f(x_i) = g_i.$ The composites
$e \circ E$, and $E \circ e$ are both the identity, and thus these
maps are invertible.

Note that if $\Gamma\to\pi$ is a surjection, it induces an inclusion
$Hom(\pi,G)\to Hom(\Gamma, G)$. In particular, if $\pi$ is a
finitely generated group, let $F_n\to\pi$ denote a fixed surjection;
endow the space of homomorphisms with the natural subspace topology
$Hom(\pi, G)\subset G^n$ or in other words make the inclusion map a
homeomorphism onto its image in the Cartesian product. Given two
surjections of this type, $F_n\to\pi$ and $F_m\to \pi$, there is a
surjection $F_{n+m}\to\pi$ which factors through both of the
original surjections, and hence gives that the topology on $Hom(\pi,
G)$ does not depend on the choice of generators. More generally, a
surjection of finitely generated discrete groups $\Gamma\to\pi$
induces a homeomorphism of $Hom(\pi,G)$ onto its image in
$Hom(\Gamma, G)$. $Hom(\pi, G)$ can be regarded as a pointed space,
for which the base-point is the constant map $f:\pi \to G$ with
$f(x) = 1_G$, the identity element in $G$.

Observe that $G$ acts by conjugation on the space of homomorphisms;
the associated orbit space $Rep(\pi, G)=Hom(\pi, G)/G$ is often
called the \textsl{representation space}. Note that the projection
onto this quotient is highly singular, and can have a complicated
structure.

The space $Hom(\pi, G)$ can in fact be realized as the fixed point
set of a $\pi$ action on the pointed mapping space \cite{Lannes}.
Consider $Map_*(\pi, G)$, the space of all maps from $\pi$ to $G$
which map $1$ to $1$. Then there is an action of $\pi$ on this space
defined as follows; for $f:\pi\to G$, $g, h\in\pi$, define
\[
gf (h) = f (hg)f(g)^{-1}
\]
from which it follows that
\[
 Hom(\pi, G) = Map_*(\pi, G)^{\pi}
\]
If $\pi$ is a finite group, then the space of homomorphisms is the
fixed point set of a smooth action on a compact manifold (namely
$G^{|\pi|-1}$), hence is itself a smooth manifold, not necessarily
connected.

The following proposition records the facts above as well as some
straightforward consequences.

\begin{prop} \label{prop:connectivity.for.free.sources}
Let $\pi$ denote a finitely generated discrete group, and $G$ a Lie
group.
\begin{enumerate}
\item
If $\Gamma$ and $\pi$ are finitely generated discrete groups with
$p:\Gamma \to \pi$ a surjective homomorphism, then the induced map
$$p^*:Hom(\pi, G) \to Hom(\Gamma ,G)$$ is a homeomorphism onto its image.
In particular if $ \Gamma = F_n$, and $\pi$ is discrete, then
$p^*:Hom(\pi, G) \to Hom(F_n,G)$ is a homeomorphism onto its image
in $G^n$. Note that in this case, $Hom(\pi, G)$ is Hausdorff.

\item If
$$1\to G \to H \to K\to 1$$
is a group extension, then
$$1\to Hom(\pi, G) \to Hom(\pi, H) \to Hom(\pi,K)$$ is an exact sequence
of sets where the base-point in $Hom(\pi, G)$ is the constant map
$c_G$ with $c_G(x) = 1_G$ for all $x$ in $\pi$.

\item If $\pi$ is isomorphic to the free product of discrete groups
$A*B$, then there is a homeomorphism
$$ Hom(\pi, G) \cong Hom(A, G) \times Hom(B, G).$$

\item If $G$ is an abelian group, then the natural map
$$Hom(\pi / [\pi, \pi], G)
\to Hom(\pi, G)$$ is a homeomorphism. Hence if $\pi/[\pi, \pi] \cong
F \oplus T$ where $F$ is free abelian and $T$ is torsion, then
$$Hom(\pi, \mathbb S^1)\cong (\mathbb S^1)^{rk~F}\times T.$$

\item If $\pi$ is a finite group, then
$Hom(\pi, G)$ is a smooth submanifold of $G^{|\pi|-1}$.
\end{enumerate}

\end{prop}

Important examples arise when considering ordered $n$--tuples of
commuting elements in a Lie group $G$; they can be identified with
$Hom(\mathbb Z^n, G)$. A key question is whether or not these spaces
are path--connected.

\begin{prop}
If every abelian subgroup of $G$ is contained in a path-connected
abelian subgroup, then the space $Hom(\mathbb Z^n, G)$ is
path-connected.

\end{prop}

\begin{proof}
Assume that every abelian subgroup of $G$ is contained in a
path-connected abelian subgroup. Let $e_1,\dots, e_n$ denote the
standard basis for $\mathbb Z^n$; given any homomorphism $f: \mathbb
Z^n \to G$, consider the images of these generators:
$x_1=f(x_1),\dots , x_n=f(e_n)$.  These elements commute, hence
$x_1,x_2,\dots ,x_n$ all lie in a path-connected abelian group
$\mathbb T$. Choose paths $p_i:[0,1]\to \mathbb T\subset G$ from
$p_i(0)=1$ to $p_i(1)=x_i$ for all $i=1,\dots , n$.
Define a homotopy $H:[0,1]\to Hom (\mathbb Z^n, G)$ by the formula
$H(t)(e_i) = p_i(t)$; by multiplying, this extends to a
well--defined homotopy in the space of homomorphisms between the
trivial homomorphism and our given homomorphism $f$.
\end{proof}

\noindent Note that if $G$ itself is abelian, then there is an
equivalence $Hom(\mathbb Z^n, G)\cong G^n$.

%The next corollary follows from the classification of compact simple
%Lie groups and thus applies to $SU(n)$, $Sp(n)$, $G_2$, $F_4$,
%$E_6$, $E_7$, or $E_8$. The conclusion of the corollary is not
%satisfied for $n$ sufficiently large with $G = SO(n)$ or $Spin(n)$.

\begin{cor}
If $G$ is equal to $U(n)$, $SU(n)$ or $Sp(n)$, then $Hom(\mathbb
Z^n, G)$ is path--connected.
\end{cor}

The conclusion of the corollary is not satisfied for $n$
sufficiently large with $G = SO(n)$ or $Spin(n)$.

Let $Map_*(B\pi, BG)$ denote the space of based maps between the
classifying spaces of $\pi$ and $G$ respectively. Its set of
path--components is denoted by $[B\pi, BG]$ and corresponds to based
homotopy classes of maps between the two classifying spaces. The
following result is well--known.

\begin{lem}\label{lem:B}
The natural map $B:Hom(\pi, G)\to Map_*(B\pi, BG)$ given by
$f\mapsto Bf$ is a continuous, and there is an induced map of sets
$$B_0: \pi_0(Hom(\pi, G))\to [B\pi, BG].$$
These maps factor through the quotient $Rep(\pi, G)=Hom(\pi, G)/G$,
and the natural map $\pi_0(Hom(\pi, G)) \to \pi_0(Rep(\pi, G))$ is a
bijection of sets is $G$ is path--connected.

\end{lem}

Note that the classifying space functor is continuous as long as the
compactly generated topology is used on the source, which is the
case here.
Clearly it makes sense to use $B_0$ to distinguish path--components
in the space of homomorphisms. More generally the following basic
question arises:
\bigskip

\noindent\textbf{Problem}:
Give conditions on $\pi$ and $G$ which imply that
$B_0:\pi_0(Hom(\pi, G)) \to [B\pi, BG]$ is a bijection.

\bigskip

The motivation of course is to understand the contribution of
representation theory to principal $G$--bundles over $B\pi$. In
special circumstances e.g. when $B\pi$ is a compact manifold, this is
geometrically significant.

\begin{exm}

Consider the case when $\pi = \pi_1(M_g)$, and $M_g$ is a compact
orientable Riemann surface of genus $g>1$. By \cite{Li}, \cite{HL},
if $G$ is a compact, connected, semisimple Lie group, there is a
bijection $\pi_0(Hom(\pi, G))\cong H^2(M_g, \pi_1(G))$. Note that in
this example, $M_g = B\pi$; given a homomorphism $f:\pi \to G$,
there is an induced map $Bf: M_g\to BG$. The bijection is obtained
by looking at the map induced in homology $H_2(M_g,\mathbb Z)\to
H_2(BG,\mathbb Z)\cong \pi_2(BG)\cong \pi_1(G)$; note that it
factors through the map $B_0$. Moreover, using the standard
cofibration
$$\mathbb \bigvee_{i=1}^{2g}~\mathbb S^1\to M_g\to \mathbb S^2,$$ it
follows that there is a bijection $[M_g,BG]\cong \pi_1(G)$, whence
the map $B_0$ must also be a bijection. In contrast, for the case of
genus $g=1$, recall that $Hom(\mathbb Z^2, U(n))$ is path connected,
even though $\pi_1(U(n))\ne 1$. Hence in this case the map $B_0$ is
trivial and so fails to be bijective.
\end{exm}

There is a class of groups built out of free abelian groups, and
which are geometrically quite significant. The definition below
appears in \cite{ACC}:

\begin{defin}\label{defin:homologically.toroidal}

A discrete group $\pi$ is said to be homologically toroidal if there
exists $H$, a finite free product of free abelian groups of finite
rank as well as a homomorphism $\phi: H\to \pi$ which induces a
split epimorphism in integral homology.

\end{defin}

The standard examples of such groups are fundamental groups
associated to arrangements of hyper-planes. More precisely, let
$\mathcal A$ denote a finite arrangement of affine hyper-planes in
$\mathbb C^n$, and let $\mathcal{M}\cong \mathbb
C^n~\backslash~\cup_{H\in\mathcal{A}}~H $ denote their complement.
Then, if this complement has contractible universal cover, $\pi =
\pi_1(\mathcal{M})$ is a homologically toroidal group. The best
known examples here are the pure braid groups $P_n$.

The next result was stated in the introduction as Theorem
\ref{thm:homologically.toroidal.reps}.

\begin{thm}
Let $\pi$ denote a homologically toroidal group; then, if $G$ is a
Lie group with path--connected maximal abelian subgroups, the map

$$B_0: \pi_0(Hom(\pi, G))\to [B\pi, BG]$$
is trivial. In particular any complex unitary representation of such
groups induces a trivial vector bundle over $B\pi$.
\end{thm}

\begin{proof}

There is a commutative diagram:

\[
\begin{CD}
  \pi_0(Hom(\pi, G)) @>{B_0}>> [B\pi, BG]      \\
  @VV{\phi^*}V    @VV{B\phi^*}V       \\
\pi_0(Hom(H, G)) @>{B_0'}>> [BH, BG]\\
\end{CD}
\]

\bigskip

Note that as $H\cong H_1*\dots *H_r$, where each $H_i$ is free
abelian of finite rank, then $Hom(H,G)\cong Hom(H_1,
G)\times\dots\times Hom(H_r,G)$; this is a product of path connected
spaces--hence it is also path--connected. Therefore the composition
$B\phi^*\cdot B_0$ is trivial. However, as was shown in \cite{ACC},
page 529, the inverse image of every point under $B\phi^*$ consists
of a single element, as this is just the Puppe exact sequence
induced by the cofibration $BH\to B\pi\to C_{B\phi}$ where the
second map is null-homotopic by our hypotheses. Hence it follows
that $B_0$ is trivial, verifying the claim.
\end{proof}

\begin{cor}
If $\mathcal{M}$ is the complement of an arrangement of hyperplanes
which is aspherical, then every unitary representation of
$\pi_1(\mathcal{M})$ induces a trivial bundle over $\mathcal{M}$.
\end{cor}

Next consider target groups with maximal abelian subgroups which are
not connected. The basic examples are the orthogonal groups $O(n)$,
the special orthogonal groups $SO(n)$ and their double covers, the
spinor groups $Spin(n)$.

If $O(n)$ is the usual real orthogonal group, then there is a split
extension

$$1\to SO(n)\to O(n)\to\mathbb Z/2\to 1.$$
The double cover of $SO(n)$ is called $Spin(n)$; it fits into a
central group extension

$$1\to \mathbb Z/2\to Spin(n)\to SO(n)\to 1.$$
When $n=3$, there is the classical extension

$$1\to \mathbb Z/2\to SU(2)\to SO(3)\to 1.$$
Hence there are exact sequences of pointed sets

$$1\to Hom(\pi, SO(n))\to Hom(\pi, O(n))\to
Hom(\pi, \mathbb Z/2)$$ and

$$1\to Hom(\pi,\mathbb Z/2) \to Hom(\pi,Spin(n))) \to
Hom(\pi,SO(n)).$$ Note that the term $Hom(\pi,\mathbb Z/2)$ is equal
to $H^1(\pi,\mathbb Z/2)$ (mod $2$ group cohomology).

As observed earlier, $O(n)$ is a semidirect product $O(n)\cong
SO(n)\times_T\mathbb Z_2$ (in fact for $n$ odd this is a direct
product). Using a section $\mathbb Z/2\to O(n)$ for the determinant
map, the following can be established:

\begin{lem}
The determinant map
$$\psi: Hom(\pi, O(n))\to H^1(\pi, \mathbb F_2)$$
is a split surjection, with each
$$\psi^{-1}(w) = Hom(\pi, O(n))_w$$
a non--empty open (and closed) subspace. This decomposes $Hom(\pi,
O(n))$ into disjoint subspaces:

$$Hom(\pi, O(n))\cong
\bigsqcup_{w\in H^1(\pi,\mathbb F_2)} Hom(\pi, O(n))_w$$ and in
particular
$$\# \pi_0(Hom(\pi, O(n))=\sum_{w\in H^1(\pi,\mathbb F_2)}
                        \#\pi_0(Hom(\pi, O(n))_w).$$
\end{lem}

\begin{cor}
If $H^1(\pi, \mathbb F_2)\ne 0$, then $Hom(\pi, O(n))$ is not
path--connected.
\end{cor}

The subspace $Hom(\pi, O(n))_w\subset Hom(\pi, O(n))$ will be
referred to as the $w$--sector of the space of homomorphisms.

Consider the Stiefel--Whitney classes $w_i(f)=f^*(w_i)\in
H^i(\pi,\mathbb F_2)$ associated to an orthogonal representation
$f:\pi \to O(n)$. These are homotopy invariants of $Bf:B\pi\to
BO(n)$, hence
they are well defined functions $[B\pi, BO(n)] \to H^i(\pi, \mathbb
F_2)$ which can, in some situations, detect path components for the
space of homomorphisms. Note that $f:\pi\to O(n)$ lands in $SO(n)$
if and only if $w_1(f)=0$, and further it will lift to $Spin(n)$ if
and only if $w_2(f)=0$. The subspace $Hom(\pi, O(n))_w$ consists of
all homomorphisms $f:\pi\to O(n)$ with first Stiefel--Whitney class
equal to $w$, and the $0$--sector is precisely $Hom(\pi, SO(n))$.

The next result makes use of the second Stiefel--Whitney class.

\begin{lem}\label{lem:w2}
Let $A$ be a finitely generated abelian group and $w$ a fixed class
in $H^1(\pi, \mathbb F_2)$. Then for $n$ sufficiently large, the
natural map
$$Hom(A, O(n))_w\to H^2(A, \mathbb F_2)$$
given by $f\mapsto w_2(f)$ is surjective.
\end{lem}
\begin{proof} The proof of this for $w=0$ (i.e.
$Hom(A, SO(n))$) is routine, and follows easily from analogues of
the arguments given in \cite{ACC}, pp. 532--533. Now for general
$w\ne 0$ observe that there is a one--dimensional representation
$\epsilon: A\to O(1)$ with first Stiefel--Whitney class equal to
$w$. Given $v\in H^2(A,\mathbb F_2)$, realize it as $w_2(f)$ for
some $f:A\to SO(m)$. Now let $f'$ denote the composition $i\cdot
(f\times\epsilon)$, where $i$ is the natural inclusion $SO(m)\times
O(1)\to O(m+1)$ given by

$$(A, t)\mapsto
\begin{pmatrix}
A & 0 \\
0 & t \\
\end{pmatrix}.
$$
A simple calculation shows that $w_1(f') = w_1(\epsilon) + w_1(f)=
w$ whereas $w_2(f') = w_2(f) + w_1(f)w_1(\epsilon) + w_2(\epsilon) =
v$.
\end{proof}

This lemma can used to prove the following result.

\begin{thm}
Let $\pi$ denote a finitely generated discrete group; then for $n$
sufficiently large and $w\in H^1(\pi, \mathbb F_2)$,

$$\#\pi_0 (Hom(\pi, O(n))_w)\ge N(\pi)$$
where $N(\pi) = |im [H^2(\pi/[\pi,\pi],\mathbb F_2) \to
H^2(\pi,\mathbb F_2)]|$.
\end{thm}
\begin{proof}

Consider the commutative diagram associated to the abelianization
map $$\rho:\pi \to \pi/[\pi,\pi]$$

\[
\begin{CD}
  \pi_0(Hom(\pi/[\pi,\pi],O(n))_w) @>{w_2}>>
H^2(\pi/[\pi,\pi],\mathbb F_2)\\
  @VV{\rho^*}V    @VV{B\rho^*}V       \\
\pi_0(Hom(\pi, O(n))_w) @>{w_2}>> H^2(\pi, \mathbb F_2)\\
\end{CD}
\]

By the previous lemma, the top arrow is a surjection for $n$
sufficiently large.
It follows that the bottom arrow surjects onto the image of
$B\rho^*: H^2(\pi/[\pi,\pi],\mathbb F_2)\to H^2(\pi,\mathbb F_2)$.
\end{proof}

\begin{cor}
Let $\pi$ be a finitely generated discrete group such that
$$H^2(\pi/[\pi,\pi],\mathbb F_2)\to H^2(\pi, \mathbb F_2)$$
is surjective; then for $n$ sufficiently large,
$$\#\pi_0(Hom(\pi, O(n)))\ge
|H^1(\pi, \mathbb F_2)||H^2(\pi,\mathbb F_2)|.$$
\end{cor}

These results apply particularly well to fundamental groups of
aspherical complements of hyperplane arrangements. Their cohomology
is torsion free and generated by $1$--dimensional classes, hence in
fact for these groups their cohomology is a quotient of the
cohomology of their abelianization.

\begin{prop}
Let ${\mathcal M}$ denote the complement of an arrangement of
hyperplanes which happens to be aspherical; if $\pi$ denotes its
fundamental group and $n$ is sufficiently large, then

$$\#\pi_0(Hom(\pi, O(n))\ge
|H^1(\pi,\mathbb F_2)||H^2(\pi, \mathbb F_2)|.$$
\end{prop}

So far, the contribution of the first two Stiefel--Whitney classes
to the number of path-components of $Hom(\pi, O(n))$ has been
analyzed. However, this is very far from being the only situation
where the space of homomorphisms fails to be path-connected. Indeed
it has been shown in \cite{Kac-Smilga} that $Hom(\mathbb Z^3,
Spin(7))$ is not path-connected, even though (using methods similar
to those outlined before, see \cite{ACC}) the map
$$\pi_0(Hom(\mathbb Z^n, Spin (m)))\to [B\mathbb Z^n,
BSpin(m)]$$ is trivial, for all $m, n\ge 1$.

\section{Examples}

One basic fact about spaces of homomorphisms is that they can change
rather drastically even when the target groups are closely related.
In particular if $K\subset G$ is a maximal compact subgroup, the set
of path-components may change when considering $Hom(\pi, K)\to
Hom(\pi, G)$. This will be illustrated with examples where $\pi$ is
allowed to be more general than a finitely generated free abelian
group.

Let $G = SL(2, \mathbb R)$ with maximal compact subgroup given by
$SO(2)$. Note that $SO(2)$ is abelian while $SL(2,\mathbb R)$
contains a free group on two letters (even as a discrete subgroup).
Thus,
there is a homeomorphism
$$Hom(H_1(\pi), SO(2)) \cong Hom(\pi, SO(2)),$$
but the natural map
$$Hom(H_1(\pi), SL(2,\mathbb R)) \to Hom(\pi, SL(2,\mathbb R))$$
is not a surjection on path-components in the case when $\pi$ is the
fundamental group of a Riemann surface of genus greater than one;
notice that there are faithful representations of surface groups in
$SL(2,\mathbb R)$ by classical work of Fricke and Klein
\cite{FrickeKlein} (this also follows from \cite{Li}).

Now let $\Gamma_g$ denote the mapping class group of a closed
orientable surface of genus $g$, and let $B_n$ denote Artin's braid
group on $n$ strands.
Consider homomorphisms into the group $PGL(2, \mathbb C)$, which has
maximal compact subgroup given by $SO(3)$,

The following result has classical origins, but does not appear to
be standard.

\begin{prop}\label{prop:SO(3)}
Assume that $\pi$ is
the mapping class group $\Gamma_g$ for $g \geq 2$, or $B_n$ for $ n
\geq 6$.
\begin{enumerate}

     \item If $G = SO(3)$, any homomorphism $\rho: \pi \to G$
factors through a cyclic group, and thus the natural map
$$Hom(H_1(\pi), SO(3)) \to Hom(\pi, SO(3))$$ is a homeomorphism.

\item The space $Hom(B_n, SO(3))$ is homeomorphic to $SO(3)$ for $n > 5$.

\item The space $Hom(\Gamma_g, SO(3))$ is a point for $g > 2$, and is
$Hom(\mathbb Z/20\mathbb Z,SO(3))$ in case $g =2$.

\end{enumerate}
\end{prop}

\begin{proof}

Recall that the first homology group of the braid group $B_n$ is
$\mathbb Z$ for $n>1$ while the first homology group of the mapping
class group $\Gamma_g$ is trivial for $g > 2$ with $H_1(\Gamma_2) =
\mathbb Z /20 \mathbb Z$. Thus, $Hom(H_1(\pi), SO(3))$ is either
$SO(3)$ in case $\pi = B_n$, a point in case $\pi = \Gamma_g$ for $g
> 2$, or $Hom(\mathbb Z/20\mathbb Z, SO(3))$ a closed subspace of
$SO(3)$ in case $\pi = \Gamma_2$. The abelianization map induces an
embedding onto a closed subspace of $Hom(\pi, SO(3))$.
To prove part (1), it suffices to see that this map is surjective.

First assume that $\pi$ is $B_n$ for $n \geq 6$. The braid groups
have the following classical presentation:

\begin{itemize}
\item $\sigma_i \cdot  \sigma_j = \sigma_j \cdot  \sigma_i$ if
$|i-j| > 2$ with $1 \leq i,j \leq n-1$ and
\item $\sigma_i \cdot  \sigma_{i+1} \cdot \sigma_i = \sigma_{i+1}\cdot \sigma_i \cdot
\sigma_{i+1}$ with $1 \leq i < i+1 \leq n-1$.
\end{itemize}

\noindent Thus all generators are conjugate. Hence if $\rho: \pi \to
SO(3)$ is a homomorphism, all elements $\rho(\sigma_i)$ are
rotations through equal angles $\theta$, with axis of rotation given
by $\mathbb X_i$.

Assume that $\theta \neq \pi, 0$. Since $\sigma_1$ commutes with
$\sigma_i$ for $ i > 2$, it follows that $\rho(\sigma_i)$ is in the
maximal torus containing $\rho(\sigma_1)$, $\mathbb T_1$. Thus the
axis of rotation $\mathbb X_1$ is equal to $\mathbb X_i$ for $i>2$,
the rotations are equal, and $\rho(\sigma_1) = \rho(\sigma_i)$ for
$i > 2$. A similar argument gives that $\rho(\sigma_2) =
\rho(\sigma_i)$ for $i > 3$. Hence if $\theta \neq \pi, 0$, then
$\rho(\sigma_1) = \rho(\sigma_i)$ for $i > 1$, and the first two
statements follow.

It remains to check the cases of $\theta = 0$, and $\theta = \pi$
with the first case being trivial. If $\theta = \pi$, then assume
that $\rho(\sigma_1)$ is rotation through the $z$-axis, and is given
by the matrix

\[
\begin{pmatrix}
-1      & \hfill 0  & \hfill 0 \\
  \hfill 0  & -1   & \hfill 0 \\
   \hfill 0  & \hfill 0   & \hfill 1 \\
\end{pmatrix}.
\] Since $\rho(\sigma_1)$ commutes with $\rho(\sigma_i)$
for $i > 2$, $\rho(\sigma_i)$ is equal to

\[
\begin{pmatrix}
a(i)      &  b(i)    & \hfill 0 \\
c(i)    & d(i)     & \hfill 0 \\
   \hfill 0  & \hfill 0   & e(i)  \\
\end{pmatrix}.
\] Since each of these matrices is in $SO(3)$, and
the square of these matrices are the identity as $\theta$ is
rotation through $\pi$, $e(i)^2 = 1 =a(i)^2 = d(i)^2$ with $c(i) = 0
= b(i)$. Thus these matrices are diagonal with entries $\pm 1$, the
image $\rho(\sigma_i)$ for $i \neq 2$ lies in an abelian group, the
$2$-torus $(\mathbb Z/ 2 \mathbb Z)^2$. By the relation $\sigma_i
\cdot \sigma_{i+1} \cdot \sigma_i = \sigma_{i+1} \cdot \sigma_i
\cdot \sigma_{i+1}$, it follows that $\rho(\sigma_1) =
\rho(\sigma_i)$ for $i \geq 3$.

A similar argument using the fact that $\sigma_i$ commutes with
$\sigma_2$ gives that $\rho(\sigma_2) = \rho(\sigma_i)$ for $i>3$.
Thus $\rho(\sigma_1) =  \rho(\sigma_j)$ for $j \geq 2$ and the
natural map $$Hom(H_1(\pi), SO(3)) \to Hom(\pi, SO(3))$$ is a
continuous bijection which must be a homeomorphism as both spaces
are compact and Hausdorff.

A similar argument works for $\Gamma_g$ with $g \geq 2$. Namely, the
mapping class group $\Gamma_g$ is generated by $\sigma_i$ for $1
\leq i \leq 2g+2$ plus one more element $\delta$ which commutes with
$\sigma_j$ for $j \neq 4$ while all the generators are conjugate (by
B. Wajnryb's presentation, see \cite{W}). Parts (2) and (3) follow
from the above.
\end{proof}

Next, consider $Sp(4,\mathbb R)$ with maximal compact subgroup equal
to $U(2)$.

\begin{thm}\label{thm:maximal.compact.in.G}
The natural map $$Hom(\Gamma_2, U(2))\to Hom(\Gamma_2, Sp(4, \mathbb
R))$$ fails to be a surjection on path-components.
\end{thm}

\begin{proof}

The method of proof is to use the results in \ref{prop:SO(3)}
together with properties of the map
$$B:Hom(\pi, G) \to Map_*(B\pi, BG)$$ which sends a homomorphism
$g$ to $Bg$ together with the induced map on the level of
cohomology.

Consider the composite
  \[
\begin{CD}
\Gamma_2  @>{f}>> Sp(4, \mathbb Z)@>{i}>> Sp(4, \mathbb R)
\end{CD}
\] where $i:Sp(4, \mathbb Z)\to Sp(4, \mathbb R)$ is the natural
inclusion. Assume that this composite lifts up to conjugacy to the
maximal compact subgroup $U(2)$ of $Sp(4,\mathbb R)$ to give a
commutative diagram
\[
\begin{CD}
\Gamma_2 @>{i\circ f}>> U(2)  \\
 @VV{1}V          @VV{}V \\
\Gamma_2 @>{i\circ f}>> Sp(4, \mathbb R).
\end{CD}
\]

There is an induced commutative diagram
\[
\begin{CD}
\Gamma_2  @>{i \circ f }>> U(2) @>{\pi}>> SO(3) \\
@VV{1}V          @VV{i}V    @VV{}V \\
\Gamma_2  @>{i \circ f }>> Sp(4, \mathbb R) @>{\pi}>> PSp(4, \mathbb
R)
\end{CD}
\]

By Proposition \ref{prop:SO(3)}, any map
\[
\begin{CD}
\Gamma_2 @>{}>> SO(3)
\end{CD}
\] factors through a cyclic group. Thus the image in mod-$2$ cohomology
in any dimension is of rank at most one.

By the results in \cite{Cohen-AJM} restricting to a semi-dihedral
subgroup of order sixteen, $SD_{16}$, and the Klein $4$-group in
$\Gamma_2$, the composite $$B\Gamma_2 \to BSp(4, \mathbb R) \to
BPSp(4, \mathbb R)
$$ satisfies the property that $w_2^3$, and $w_3^2$ are linearly
independent in $H^6(B\Gamma_2;\mathbb F_2)$. Thus there does not
exist a representation $\Gamma_2 \to U(2)$ which is homotopic to
$i\circ f$ in $Hom(\Gamma_2,Sp(4, \mathbb R))$. The theorem follows.
\end{proof}

\section{Some Homological Computations}
Let $A$ denote a discrete group, and $G$ a Lie group. Consider the
space $Hom(\mathbb Z\times A, G)$. Using the surjection $\mathbb Z *
A\to \mathbb Z\times A$, the space $Hom(\mathbb Z\times A, G)$ can
be regarded as a subspace given by
$$Hom (\mathbb Z\times A, G)\subset G\times Hom (A,G).$$

\begin{prop}

Let $G$ denote a non--abelian Lie group and $A$ a finitely generated
discrete group such that $A/[A,A]\ne\{1\}$. Then the natural
projection map $pr_1: G\times Hom (A, G)\to G$ restricts to a map
$$\pi : G\times Hom(A,G)-Hom(\mathbb Z\times A, G)\to G$$
with the following properties:
\begin{itemize}
\item the image of $\pi$ is $G-Z(G)$, where $Z(G)\subset G$ is the
center of $G$

\item if $x_0\in G-Z(G)$, then

$$\pi^{-1}(x_0)\cong Hom(A,G)- Hom(A, Z_G(x_0))$$
where $Z_G(x_0)$ is the centralizer of $x_0$ in $G$.
\end{itemize}
\end{prop}
\begin{proof}
Note that if $x_0\in Z(G)$, then for any $f\in Hom (A,G)$, $im
f\subset Z_G(x_0)$, hence $x_0\notin im \pi$. Now if $x_0\notin
Z(G)$, then for all $m\ge 0$ there exists an element $u_m\in G$ of
order $m$, such that $x_0$ does not commute with $u_m$. By
the hypotheses, there exists an $m>0$ and a surjection $A\to \mathbb
Z/m$, hence a homomorphism $A\to G$ which factors through the
subgroup generated by $u_m$, such that the image of $f$ is not
contained in $Z_G(x_0)$. Hence the pair $(x_0, f)$ defines an
element in the domain of $\pi$, and so the first assertion has been
verified.

For the second part, note that if $x_0\in G-Z(G)$, then

$$\pi^{-1}(x_0) = \{(x_0,f)~|~ im f\not\subset Z_G(x_0)\}$$
which is clearly homeomorphic to the complement $Hom (A, G)-Hom (A,
Z_G(x_0))$.
\end{proof}

This proposition can be applied to an important example.

\begin{exm}
Let $A=F_r$, the free group on $r$ generators. In this case
there is a projection
$$\pi : G^{r+1}-Hom (\mathbb Z\times F_r, G)\to G-Z(G)$$
with $\pi^{-1}(x_0) \cong G^r- [Z_G(x_0)]^r$.
\end{exm}

In general this map can be quite complicated; the different strata
of the base, corresponding to elements with different centralizers,
requires some rather intricate gluing to understand the homology of
the domain. However
by choosing $G=SU(2)$, the situation simplifies. The main point is
that there are exactly two kinds of elements in $SU(2)$: the
singular ones, which are central (namely $Z(SU(2))=\{\pm I\}$) and
the regular elements $x_0$, which have a maximal torus as
centralizer, i.e. $Z_G(x_0)=\mathbb S^1$. Under these conditions,
the map $\pi$ is locally trivial.

\begin{prop}
For $A$ as above, $\pi$ is a locally trivial bundle

$$\pi : SU(2)\times Hom (A, SU(2))-Hom (\mathbb Z\times A, SU(2))
\to SU(2)-\{\pm I\}$$ with fiber $F\cong Hom (A, SU(2))-Hom
(A,\mathbb S^1)$.
\end{prop}
Note that the base of this bundle has the homotopy type of $\mathbb
S^2$.

\begin{exm}
Let $F_r$ denote the free group on $r$ generators. Then
$$\pi : SU(2)^{r+1} - Hom (\mathbb Z\times F_r, SU(2))
\to SU(2) -\{\pm I\}$$ is a locally trivial bundle with fiber
$SU(2)^r-\mathbb T^r$, where $\mathbb T^r$ is the torus of rank $r$.
Taking $r=1$ yields a fibration up to homotopy
$$SU(2) - \mathbb S^1 \to SU(2)\times SU(2)
- Hom (\mathbb Z\oplus\mathbb Z, G)\to \mathbb S^2$$ for which the
total space is
homotopy equivalent to $SO(3)$.

A geometric argument is required to justify this last statement
as well as to supply details necessary for homology calculations.
Consider the commutator map
$$\partial: SU(2)\times SU(2)\to SU(2).$$
Given the homeomorphism $Hom(\Z *\Z, SU(2))\cong SU(2)\times SU(2)$,
$Hom (\Z\oplus\Z, SU(2))$ can be identified with the subspace
$\partial^{-1}(I)$. Let ${\mathcal S}\subset Hom (\Z * \Z, SU(2))$
denote the \emph{reducible} representations. These are those pairs
of representations which can be simultaneously conjugated into the
maximal torus, in this case the diagonal matrices.
Consider the following basic fact from Akbulut-McCarthy, \cite{AM}
III.2.1, pp. 53--54:

\begin{prop}
$\partial$ is a surjective map, and the set of critical points for
$\partial$ is precisely equal to ${\mathcal S}$.
\end{prop}

Now note that in fact ${\mathcal S}=\partial^{-1}(I)$, as
\textsl{any} pair of commuting matrices in $SU(2)$ can be
simultaneously diagonalized. This is exactly the subspace of
singular points and hence it follows that
\begin{prop}
$M=SU(2)\times SU(2) \setminus{\mathcal S}$ is an open smooth
6--dimensional manifold.
\end{prop}

Next
observe that $SO(3) = SU(2)/\{\pm I\}$ acts via conjugation on $M$,
and as observed by Akbulut-McCarthy, this action is free, and the
orbit space $M/SO(3)$ is a $3$--dimensional manifold. The next
result in \cite{AM}, (VI.1.1) determines the homotopy type of this
complement.

\begin{prop}
Let $R_- = \partial^{-1}(-I)$, then

$$\partial_| : SU(2)\times SU(2) \setminus Hom (\Z\oplus\Z, SU(2))
\to SU(2) \setminus \{I\}$$ is a locally trivial bundle, with fiber
$R_-$.
\end{prop}
Note that for general position reasons, the total space must be
connected. Hence
it follows that $R_-$ is connected. Now note that the free $SO(3)$
action restricts to a free action on $R_-$, with orbit space a
single point. In other words $R_-\cong SO(3)$ (see \cite{AM}, pp
154).

\begin{cor} \label{cor:SO(3)}
There is a homotopy equivalence
$$SU(2)\times SU(2) \setminus Hom (\Z\oplus\Z, SU(2))
\simeq SO(3).$$
\end{cor}

The next ingredient is crucial for subsequent calculations:

\begin{lem}
The map in homology induced by the inclusion

$$i_*: H_3(R_-,\Z)\to H_3(SU(2)\times SU(2),\Z)$$
is zero.
\end{lem}
\begin{proof}
As in \cite{AM}, $R_- = \{(A,B)\in SU(2)\times SU(2) ~|~
ABA^{-1}B^{-1}=-I\}$. This condition implies that $tr(A)=tr(B)=0$.
Now, identify
$$\mathbb S^2 =\{A\in SU(2) ~|~ tr(A)=0\},$$

$$R_-\subset\mathbb S^2\times\mathbb S^2\subset\mathbb
SU(2)\times SU(2)$$ and the result follows.
\end{proof}

The space of homomorphisms $Hom(\Z\oplus\Z, SU(2))$ is an analytic
variety (see \cite{Goldman}, page 568) hence locally contractible.
Therefore Poincar\'e--Lefschetz Duality (see \cite{Dold}, pp. 292)
can be applied to the pair 
$(SU(2)\times SU(2), Hom (\Z\oplus\Z,
SU(2)))$, to conclude that there is an isomorphism for all $i\ge
0$:

$$H^i(Hom (\Z\oplus\Z, SU(2)),\Z)\cong
H_{6-i}(SU(2)\times SU(2), SU(2)\times SU(2) \setminus Hom
(\Z\oplus\Z, SU(2)), \Z).$$
The complement can be identified up to homotopy with $R_-\simeq
SO(3)$. From the lemma above and the long exact sequence in homology
for the pair $(SU(2)\times SU(2), R_-)$, it follows that

\[ H^i(Hom (\Z\oplus\Z, SU(2)), \Z)
\cong\left\{\begin{array}{r@{\quad\hbox{if}\quad}l}
\Z & i=0\\
0  & i=1\\
\Z & i=2\\
\Z\oplus\Z & i=3\\
\Z/2\Z & i=4\\
0 & i\ge 5\end{array}\right.\]

\end{exm}

\begin{exm}\label{P_3}
Consider now the case $A=F_2$, the free group on two generators.
Note that the pure braid group $P_3\cong \mathbb Z\times F_2$. There
is a fibration sequence

$$SU(2)\times SU(2) - \mathbb T^2
\to (SU(2))^3 - Hom (\mathbb Z\times F_2, SU(2)) \to \mathbb S^2$$
which can be used to compute the cohomology of the complement.

\begin{prop}
The Serre spectral sequence for the fibration above collapses at
$E_2$ and the cohomology of the total space is torsion--free with
Poincar\'e series
$$(1+t^2)(1+3t^3+2t^4)=1 + t^2 + 3t^3 + 2t^4 + 3 t^5
+ 2t^6.$$
\end{prop}

Using duality,
it follows that
\begin{thm}
$H^*(Hom (\mathbb Z\times F_2, SU(2))$ is torsion free, with
Poincar\'e series
$$1 + 2t^2 + 6t^3 + 2t^4 + t^5 + 2t^6.$$
\end{thm}

\end{exm}
The case for the complement of $Hom(\mathbb Z\times F_r, SU(2))$ can
be handled similarly. Let $p_r(t)$ denote the Poincar\'e series for
the rational cohomology of $SU(2)^r-\mathbb T^r$, then

\begin{thm}
The spectral sequence for the fibration
$$SU(2)^r-\mathbb T^r\to SU(2)^{r+1}-Hom (\mathbb Z\times F_r,
SU(2))\to\mathbb S^2$$ collapses at $E_2$ and the total space has
torsion--free homology and Poincar\'e series equal to
$(1+t^2)p_r(t)$.
\end{thm}

Note: $p_r(t)$ is not hard to compute. Using duality it should be
possible to obtain the Betti numbers for $Hom(\mathbb Z \times F_r,
SU(2))$.

\section{The Variety of Commuting Triples in $SU(2)$}
\label{sec:The Variety of Commuting Triples in $SU(2)$}

This section will focus on the case of commuting triples in $SU(2)$.
As before it seems best to first consider its complement. Let $X[m]$
denote the space of non--commuting $m$--tuples,

$$X[m] = (SU(2))^m - Hom (\mathbb Z^m, SU(2)).$$
Notice that the symmetric group $\Sigma_m$ acts naturally on both
the space of commuting and non--commuting $m$--tuples. Moreover it
has been shown that $X[2]$ is homotopy equivalent to $SO(3)$.

Let $\{i,j, k\} = \{ 1,2, 3\}$ and denote

$$X_{j,k} = \{ (z_1, z_2, z_3)\in (SU(2))^3 ~|~ z_jz_k\ne z_kz_j\}.$$

Note that $X[3] = X_{1,2}\cup X_{1,3}\cup X_{2,3}$. The cohomology
of $X[3]$ will be computed by using this decomposition. Note that 
each of the three pieces is homeomorphic to
$X[2]\times SU(2)$, and that $\Sigma_3$ permutes these three
pieces. The triple intersection consists of all pairwise
non--commuting triples, and is invariant under the action of the
symmetric group.
Using results in the previous section, it follows that

$$X_{i,j}\cong X[2]\times SU(2)\simeq SO(3)\times SU(2).$$
The cohomology of the intersections is computed next.
Let $X_{123}= X_{1,2}\cap X_{1,3}\cap X_{2,3}$.

\begin{prop}
The projection onto the $j,k$ coordinates
 $\pi_{jk}: X_{j,k}\to X[2]$ restricts to
fibrations
$$F_i\to X_{j,k}\cap X_{i,k}\to X[2]$$
$$E_i\to X_{123}\to X[2]$$
where the fibers can be described as
$$F_i\cong SU(2) - Z_G(z_k)\cong SU(2) - \mathbb S^1
\simeq \mathbb S^1$$
$$E_i\cong SU(2) - [Z_G(z_k)\cup Z_G(z_j)]
\cong SU(2) - (\mathbb S^1\cup \mathbb S^1)\simeq \mathbb S^1\bigvee
\mathbb S^1\bigvee\mathbb S^1.$$ The map $X[2]\to X_{123}$ given by
$\psi (z_j,z_k)_j=z_j, \psi (z_j,z_k)_k=z_k, \psi (z_j,z_k)_i =
z_jz_k$ defines a section for all three bundles.
\end{prop}
\begin{proof}
Given a triple $z_i,z_j,z_k$ of pairwise non--commuting elements,
its projection onto the $j,k$ coordinates will be all non--commuting
pairs. Further the fiber over any such pair $(z_j,z_k)$ can be
described as elements $z_i\in SU(2)$ such that $z_i$ is not in the
centralizer of $z_j$, nor in the centralizer of $z_k$. Now given
that $z_j, z_k$ do not commute, they are both regular elements, and
thus their centralizers are maximal tori (in this case isomorphic to
a circle $\mathbb S^1$) which intersect exactly at the center of
$SU(2)$, namely $\{\pm I\}$. Hence the fiber $E_i$ is the complement
in $SU(2)$ of a union of two unknotted circles intersecting at two
points. It is easy to see that this has the homotopy type of a wedge
of three circles. Similarly the fiber $F_i$ can be identified with
the elements $z_i\in SU(2) - Z_G(z_k)\simeq \mathbb S^1$.

Note that $\pi_{j,k}$ restricted to $X_{123}$ and $X_{j,k}\cap
X_{i,k}$ respectively still maps onto $X[2]$. Furthermore the fibers
are all homeomorphic and the bundles are locally trivial--hence they
are fibrations.

Finally note that if $z_j$ and $z_k$ do not commute, then $z_j$
cannot commute with $z_jz_k$ and $z_k$ cannot commute with $z_jz_k$;
hence $\psi (z_j,z_k)$ defines a non--commuting triple and hence a
section for all three bundles.
\end{proof}

In the fibrations above, the inclusions
$$X_{123}\subset X_{j,k}\cap X_{i,k}\subset X_{j,k}$$
cover the identity map on $X[2]$ and induce the natural inclusions
on the fibers:
$$SU(2)-[\mathbb S^1\cup\mathbb S^1]\subset
SU(2) - \mathbb S^1 \subset SU(2).$$ It will be necessary to
understand the action of $\pi_1(X[2])\cong \mathbb Z/2$ on the
homology of the fibers.

The projection
$$p_1:X_{j,k}\cap X_{i,k}\to SU(2) - \{\pm I\}$$
is also a locally trivial bundle, with fibers $SU(2)-Z_G(z_k)\times
SU(2) - Z_G(z_k)\simeq \mathbb S^1\times\mathbb S^1$. From this it
follows that $\pi_1(X_{j,k}\cap X_{i,k})$ must be a quotient of
$\mathbb Z\oplus\mathbb Z$, hence abelian. Now, going back to the
bundle over $X[2]$, this implies that the action of $\pi_1(X[2])$ on
the fiber must be trivial in homology. Using the natural fibre--wise
inclusions this can be used to verify triviality of the action for
the bundle with total space $X_{123}$.

The associated Serre spectral sequences in cohomology for these
fibrations will have an untwisted $E_2$--term, and will in fact
collapse from the start, yielding\footnote{In fact, with a bit of
additional work it is possible to show that 
each total space is homotopy
equivalent to the corresponding product bundle.}

\begin{prop}
There are additive integral cohomology equivalences
$$X_{j,k} \sim SO(3)\times SU(2),
~~X_{j,k}\cap X_{i,k} \sim SO(3)\times \mathbb S^1, ~~X_{123}\sim
SO(3)\times [\mathbb S^1\bigvee\mathbb S^1\bigvee \mathbb S^1].$$
\end{prop}

\begin{proof}
The first is already known. For the other two note that the fibers
have torsion--free (co)homology, all concentrated in degree one. The
existence of a section implies collapse of the associated spectral
sequence (untwisted) at $E_2$. An explicit computation and
comparison with the spectral sequence mod 2 shows that there are no
extension problems and the result follows.
\end{proof}

Using the inclusions $\rho_i:E_i\to F_i$ it is possible to find a
basis $\{e_1, e_2, e_3\}$ for $H^*(E_i, \mathbb Z)$ such that if
$u\in H^1(F_i,\mathbb Z)$ is the canonical generator, then
$\rho_i^*(u) = e_i$. As a consequence of this analysis it follows
that the alternating sum of the maps induced by inclusions in
cohomology induces an isomorphism for $k=1$:
$$H^1(X_{1,2}\cap X_{1,3},\mathbb Z)\oplus
H^1(X_{1,2}\cap X_{2,3},\mathbb Z) \oplus H^1(X_{1,3}\cap
X_{2,3},\mathbb Z) \cong H^1(X_{123}, \mathbb Z).$$

Consider the Mayer-Vietoris spectral sequence in integral cohomology
associated to the covering $X = X_{1,2} \cup X_{1,3}\cup X_{2,3}$.
It will have as its $E_1$ term

\[ E_1^{p,q}
\cong\left\{\begin{array}{r@{\quad\hbox{if}\quad}l}
H^q(X_{1,2},\mathbb Z)\oplus H^q(X_{1,3},\mathbb Z)
\oplus H^q(X_{2,3},\mathbb Z) & p=0\\
H^q(X_{1,2}\cap X_{1,3},\mathbb Z) \oplus H^q(X_{1,2}\cap
X_{2,3},\mathbb Z)\oplus
H^q(X_{1,3}\cap X_{2,3},\mathbb Z) & p=1\\
H^q(X_{123},\mathbb Z) & p=2\\
0 & p\ge 3\end{array}\right.\]

The differentials on the $E_1$--page are determined as alternating
sums of the inclusion maps. From the previous observation
$d_1^{1,q}: E_1^{1,q}\to E_1^{2,q}$ is an epimorphism for all
$q\ge 0$. Hence there are no higher differentials in the spectral
sequence, which completely collapses at $E_1$, to yield

\begin{thm}\label{thm:commuting.triples}

$$H^*(X[3],\mathbb Z)=H^*((SU(2))^3 - Hom (\mathbb Z^3, SU(2)),\mathbb Z) \cong H^*([\bigvee_{i=1}^3 SU(2)]\times SO(3),
\mathbb Z).$$
\end{thm}

Once again using the fact that the spaces of homomorphisms under
consideration are locally contractible, and applying
Poincar\'e--Lefschetz duality to the pair 
$$(SU(2))^3, Hom (\mathbb
Z^3, SU(2))$$ 
yields (see \cite{Dold}, page 292) that for all $i\ge
0$ there is an isomorphism

$$H^i(Hom (\mathbb Z^3, SU(2)), \mathbb Z)\cong H_{9-i}((SU(2))^3,
X[3]), \mathbb Z)$$ where as before $X[3]$ is the space of
non--commuting triples, whose (co)homology has just been computed.
The basic ingredient needed is to understand the map induced in
homology by the inclusion $X[3]\subset SU(2)\times SU(2) \times
SU(2)$. This is determined by the inclusions $X_{j,k}= X[2]\times
SU(2)\subset (SU(2))^3$ and this in turn by the inclusion
$X[2]\subset SU(2)\times SU(2)$. However it has already been shown
that up to homotopy there is a factorization
$$X[2]\subset \mathbb S^2\times\mathbb S^2\to SU(2)\times SU(2)$$
where $\mathbb S^2$ is identified with matrices of trace zero. Hence
$X[2]\subset SU(2)\times SU(2)$ is trivial in reduced homology, and
taking the product with $SU(2)$ it follows that the map in reduced
homology induced by $X_{j,k}\subset (SU(2))^3$ is only non--trivial
on the three dimensional class corresponding to the sphere, which is
mapped via the identity. Putting these together, and observing that
the sphere appearing in the $X_{j,k}$'s rotates around the three
spheres in the product leads to

\begin{lem}
The map in reduced integral homology induced by the inclusion
$X[3]\subset (SU(2))^3$ is a surjection in dimension equal to
three, and zero everywhere else.
\end{lem}

\noindent Applied to the long exact sequence in homology for
the pair $((SU(2))^3, X[3])$, this yields the following computation:

\begin{thm}

The integral cohomology of the space of (ordered) commuting triples
in $SU(2)$ is given by

\[ H^i(Hom (\Z\oplus\Z\oplus\Z, SU(2)), \Z)
\cong\left\{\begin{array}{r@{\quad\hbox{if}\quad}l}
\Z & i=0\\
0  & i=1\\
\Z\oplus\Z\oplus\Z & i=2\\
\Z\oplus\Z \oplus\Z & i=3\\
\Z/2\Z\oplus\Z/2\Z\oplus\Z/2\Z & i=4\\
\Z & i=5\\
0 & i=6\\
\Z/2\Z & i= 7
\\0  & i\ge 8\end{array}\right.\]
\end{thm}
See \cite{BFM} for a different perspective on commuting elements in
compact Lie groups.

\section{Degenerate Subspaces and Stable Structure for Commuting Elements}

Natural
subspaces of $Hom(\mathbb Z^{n},G)$ arise from the fat wedge
filtration of the product $G^n$ where the base-point of $G$ is $1_G$
with $F_jG^n$ given by the subspace of $G^n$ with at least $j$
coordinates equal to $1_G$. Define subspaces of $Hom(\mathbb
Z^{n},G)$ by the formula
$$S_{n}(j,G) = Hom(\mathbb Z^{n},G)\cap F_jG^n$$
and let $S_n(G)=S_n(1,G)$. There is an induced filtration
$$S_n(n,G)\subset S_n(n-1,G)\subset\dots\subset S_n(1,G)
\subset Hom(\mathbb Z^n, G).$$
The key property required is the following.

\begin{defin}\label{defin:cofibrantly commuting}
A Lie group $G$ is said to
have cofibrantly commuting elements if the
natural inclusions
$$I_j: S_{n}(j,G)\to S_{n}(j-1,G)$$ are cofibrations with cofibre
$S_n(j-1,G)/S_n(j,G)$ for all $n$ and $j$ for which both spaces are
non-empty.
\end{defin}

A result due to Steenrod \cite{St} is that if the
pair $(X,A)$ is an NDR-pair, then the map $A \to X$ is a cofibration
with cofibre $X/A$. In section \ref{sec:Cofibrations},
it will be verified that if $G$ is
a closed subgroup of $GL(m, \mathbb C)$, then the pairs $(S_n(j-1,
G), S_n(j,G))$ are all NDR pairs and hence that these groups $G$
\textsl{have cofibrantly commuting elements}.
Given this technical condition,
it will follow that the spaces $Hom(\mathbb Z^n,G)$ naturally split
after a single suspension, where the summands which appear are the
suspensions of the successive quotients $S_n(j-1,G)/S_n(j,G)$
arising from the filtration. These in turn can be expressed using
quotients of the form $Hom(\mathbb Z^q, G)/S_q(G)$.

The decompositions of $\Sigma(Hom(\mathbb Z^n,G))$ will make use of
the following lemma.\footnote{Recall that if $X$ and $Y$ are pointed
topological spaces, then $X\bigvee Y$ denotes their one point union
in the Cartesian product $X\times Y$, while $\Sigma X$ denotes the
suspension of $X$.}

\begin{lem}\label{lem:split.cofibrations}
Let
\[
\begin{CD}
A @>{i}>> B @>{p} >> C
\end{CD}
\]
be a cofibration of path-connected spaces of the homotopy type of a
CW-complex that is split after one suspension, i.e. assume given a
map
$$r: \Sigma(B) \to \Sigma(A)$$ such that $r\circ \Sigma(i)$ is homotopic
to $1_{\Sigma(A)}$. Then the map
\[
\begin{CD}
\Sigma (B) @>{r \bigvee \Sigma(p)}>> \Sigma(A) \bigvee \Sigma(C)
\end{CD}
\] is a homotopy equivalence.

\end{lem}

\begin{proof}
Consider the morphism of cofibrations given by
\[
\begin{CD}
\Sigma(A)   @>{\Sigma(i) }>> \Sigma(B)  @>{\Sigma(p) }>> \Sigma(C)  \\
@VV{1}V                @VV{r \bigvee \Sigma(p)}V @VV{1}V
\\
\Sigma(A)   @>{\Sigma(i) }>> \Sigma(A) \bigvee \Sigma(C)
@>{project}>> \Sigma(C)
\end{CD}
\] Thus $r \bigvee \Sigma(p)$ is a homology isomorphism. This suffices.

\end{proof}

Assume that  $$i: A \to B$$ is a cofibration. The crux of the matter
is to give a choice of map $r:\Sigma(B) \to \Sigma(A)$ as given in
the lemma. An important example which will be made use of is the
inclusion $X\bigvee Y\to X\times Y$, where $X$ and $Y$ are pointed
$CW$--complexes. In this case the splitting map
$$\pi: \Sigma (X \times Y) \to \Sigma (X) \bigvee \Sigma(Y)$$ is given
by the following composite
\[
\begin{CD}
\Sigma (X \times Y) @>{pinch}>> \Sigma (X \times Y) \bigvee \Sigma
(X \times Y) @>{\Sigma (p_X) \bigvee \Sigma (p_Y) }>> \Sigma (X)
\bigvee
\Sigma(Y)\\
\end{CD}
\] where $p_X: X \times Y \to X$ is the natural projection map.

Consider the inclusion $\bigvee^nG\to Hom(\mathbb Z, G)$ and its
mapping cone, denoted $A_n(G)$. By construction this will be a
cofibration sequence. The composite
\[
\begin{CD}
\bigvee^n G @>{}>> Hom(\mathbb Z^n,G) @>{}>>  G^n,
\end{CD}
\] factors the inclusion of the wedge $\bigvee^n G$
in the product $G^n$. This map is split after a single suspension;
hence it follows that after one suspension, the natural cofibration
$\bigvee^n G \to Hom(\mathbb Z^n,G)  \to A_n(G)$ is also split. From
this it follows that if $G$ is any Lie group, then there is a
homotopy equivalence
$$\Sigma ( Hom(\mathbb Z^n,G))\simeq
\Sigma ( \bigvee^n  G ) \bigvee \Sigma(A_n(G)).$$ However, the
geometric identification of $A_n(G)$ has not been made precise, even
for $n=2$. This is precisely what is desired for certain important
examples, assuming the hypothesis that $G$
has cofibrantly commuting elements.

In what follows consideration will be given to the successive
quotients $K(n,q,G)=S_n(n-q,G)/S_n(n-q+1)$, and it will be important
to have a precise description of them.

\begin{thm} \label{thm:homeomorphisms.and.subquotients}
The quotient $$K(n,q,G) = S_{n}(n-q,G)/S_{n}(n-q+1,G)$$ is
homeomorphic to the subspace of the product $$\prod^{\binom n q}
Hom(\mathbb Z^q,G)/S_{q}(G)$$ given by the wedge
$$\bigvee^{\binom n q} Hom(\mathbb Z^q,G)/S_{q}(G).$$ Thus the
suspension $\Sigma K(n,q,G)$ is a retract of $\Sigma(\prod^{\binom n
q} Hom(\mathbb Z^q,G)/S_{q}(G))$.
\end{thm}

\begin{proof}
Given an ordered $k$ element subset $I = (n_1, n_2, \cdots, n_k)$ of
$\{1,2, \cdots,n\}$ with $$ 1 \leq n_1 < n_2 < \cdots < n_k \leq
n,$$ there are associated projection maps
$$P_I: Hom(\mathbb Z^n,G) \to Hom(\mathbb Z^k,G).$$
There are embeddings
$$S_I: Hom(\mathbb Z^k,G) \to Hom(\mathbb Z^n,G)$$
by inserting the identity $1_G$ in the coordinates not equal to
$n_1, n_2, \cdots, n_k$. The composite
\[
\begin{CD}
Hom(\mathbb Z^k,G) @>{S_I}>> Hom(\mathbb Z^n,G)@>{P_I}>> Hom(\mathbb
Z^k,G)
\end{CD}
\] is the identity. Thus if $k < n$, $Hom(\mathbb Z^k,G)$ is a
retract of $Hom(\mathbb Z^n,G)$. Notice that an analogous statement
is satisfied for any simplicial space.

These projection maps can be assembled to yield a map
$$Hom(\mathbb Z^n, G)\to \prod^{\binom n q} Hom(\mathbb Z^q,
G)/S_q(G).$$ Note that the subspace $S_n(n-q+1,G)$ will map to a
point, and that $S_n(n-q,G)$ maps onto the wedge
$$\bigvee^{\binom n q} Hom(\mathbb
Z^q,G)/S_{q}(G)$$ as a subspace of the product. It is easy to see
that this induces a continuous, closed bijection
$$K(n,q,G)\to \bigvee^{\binom n q} Hom(\mathbb
Z^q,G)/S_{q}(G)$$ and the result follows.

\end{proof}

In order to proceed any further, the aforementioned technical
condition on the target group $G$ must be assumed, i.e. that
$G$
has \textsl{cofibrantly commuting elements}.

\begin{thm}
Let $G$ be a Lie group which has cofibrantly
commuting elements,
then the inclusions
$$S_n(q,G)\to S_n(q-1, G)$$
split after a single suspension.
\end{thm}

\begin{proof}
It suffices to show that
$$S_n(n-q,G)\to Hom(\mathbb Z^n,G)$$
splits after a single suspension. Induction on $q$ will be used. For
$q=1$, $S_n(n-1,G)=\bigvee^n G$ and the splitting exists by using
the inclusion into the product $\prod^n G$. Now assume it for $q-1$
and prove it for $q$.

There is a commutative diagram

\[
\begin{CD}
S_n(n-q+1,G) @>{1}>> S_n(n-q+1,G)   \\
@VV{}V                @VV{}V          \\
S_n(n-q,G) @>>> Hom(\mathbb Z^n, G) \\
@VV{}V                @VV{}V          \\
K(n,q,G) @>{\gamma}>>  Hom(\mathbb Z^n, G)/S_n(n-q+1,G)
\end{CD}
\] By the inductive hypotheses, the vertical columns split after
suspending once. Here the fact that the two columns are cofibration
sequences, arising from NDR pairs is being used. It is only
necessary to show that the bottom map splits after a single
suspension. As before, a map can be constructed

$$Hom(\mathbb Z^n, G)\to \prod^{\binom n q} Hom(\mathbb
Z^q,G)/S_{q}(G)$$ using the projections. Now note that once again
$S_n(n-q+1,G)$ maps to a single point, hence it is well defined on
the quotient. The composition

$$\Sigma Hom(\mathbb Z^n, G)/S_n(n-q+1,G)\to
\Sigma \prod^{\binom n q} Hom(\mathbb Z^q, G)/S_q(G) \to\Sigma
K(n,q,G)$$ can be used to split our map, hence completing the proof.
\end{proof}

\begin{thm}
If $G$ is a Lie group which
has \textsl{cofibrantly commuting elements}, then there are
homotopy equivalences
$$\Sigma ( Hom(\mathbb Z^n,G))\simeq \bigvee_{1 \leq k \leq n}\Sigma (\bigvee^{\binom n k}
Hom(\mathbb Z^k,G)/ S_k(G)).$$
\end{thm}
\begin{proof}
There is a filtration of $Hom(\mathbb Z^n, G)$:

$$S_n(n, G)\subset S_n(n-1,G)\subset\dots\subset S_n(1,G)\subset
Hom(\mathbb Z^n, G)$$ where each successive inclusion is a
cofibration which is split after a single suspension. By induction
there is a splitting
$$\Sigma Hom(\mathbb Z^n, G)\simeq \bigvee_{1\le q\le n}\Sigma
K(n,q,G)$$ whence the result follows.
\end{proof}

\begin{rem}
The proof given here in the special setting of $\Sigma(Hom(\mathbb
Z^n,G))$ applies verbatim to the more general setting of a
simplicial space satisfying a certain cofibration condition. Let
$X_*$ denote any simplicial space with $n$-th space given by $X_n$.
Define
$$s_q(X_n) = \cup_{0 \leq j_1 \leq  j_2 \leq
\cdots \leq  j_q \leq n} s_{j_1} s_{j_2}\cdots s_{j_q}(X_{n-q}).$$

\begin{thm} \label{thm:stable.decompositions.for.simplicial.spaces}
If $X_*$ is a simplicial space such that
the $$s_q(X_n) \to s_{q-1}(X_n)$$ are
cofibrations for all $q$,
then there are homotopy
equivalences
$$\Sigma (X_n)\simeq \bigvee_{1 \leq k \leq n}
\Sigma (\bigvee^{\binom n k} X_k/s_k(X_k)).$$
\end{thm}

For the sake of concreteness, a proof in our special setting has
been provided, while the general case is left as an exercise. This
result is really just an observation, but a possibly useful one. For
our purposes the main content is checking that a Lie group $G$
has cofibrantly commuting elements.
\end{rem}

Examples concerning these splittings for $Hom(\mathbb Z^n,SU(2))$
are given next: assume
throughout that $SU(2)$  has cofibrantly
commuting elements.
The next examples follow from the stable decompositions above in
conjunction with the Mayer-Vietoris spectral sequence computations
given in the previous section.

\begin{exm}
$Hom( \mathbb Z^2, SU(2))$:

The assumption is being made that for $G=SU(2)$, the inclusion
$$ G\bigvee G \to Hom(\mathbb Z^2,G)$$ is
a cofibration sequence with cofiber the natural quotient; then there
is a homotopy equivalence $A_2(G)\simeq Hom(\mathbb Z^2,G)/G\bigvee
G$.

Since the map $G\bigvee G \to Hom(\mathbb Z^2,G) $ is an embedding,
it follows that $$ G\times G - Hom(\mathbb Z^2,G) \subset G\times G
- G\bigvee G
$$ is an embedding. Recall that
there is an embedding $SO(3) \to G\times G - Hom(\mathbb Z^2,G)$
which is a homotopy equivalence, and that the image of $SO(3)$ has a
neighborhood which is a product bundle $SO(3) \times \mathbb R^3$.
It follows that $A_2(SU(2))$ is homotopy equivalent to the
Spanier-Whitehead dual $\mathbb S^6 - SO(3)$.

\begin{prop}
If $G = SU(2)$, there is a homotopy equivalence
$$\Sigma ( Hom( \mathbb Z^2,G))\simeq
\Sigma (G\bigvee G) \bigvee \Sigma(\mathbb S^6 - SO(3)).$$
\end{prop}
From this a cohomology calculation follows:
$$H^i(Hom(\mathbb Z^2, SU(2)),\mathbb Z) = H^i(\bigvee^2 SU(2)\bigvee (\mathbb S^6
   -SO(3)),\mathbb Z).$$
Note that this agrees with the previous computation.
\end{exm}

\begin{exm} $Hom(\mathbb Z^3, SU(2))$:

After a single suspension, $Hom(\mathbb Z^3,SU(2))$ is homotopy
equivalent to $$\bigvee^3 SU(2) \bigvee
S_3(1,SU(2))/S_3(2,SU(2))\bigvee
  Hom(\mathbb Z^3,SU(2))/S_3(1,SU(2))$$
Now
$$S_3(1,SU(2))/S_{3}(2,SU(2)) = \bigvee^3 Hom(\mathbb
Z^2,SU(2))/S_{2}(1,SU(2))\simeq \bigvee^3 [\mathbb S^6 - SO(3)].$$
By our homological computations, $Hom(\mathbb
Z^3,SU(2))/S_3(SU(2))$ is homotopy equivalent to $SU(2) \wedge
(S^6 -SO(3))$. Hence it follows that

$$\Sigma Hom(\mathbb Z^3,SU(2))\simeq
\Sigma\bigvee^3 SU(2)\bigvee [\Sigma\bigvee^3 (\mathbb S^6 -
SO(3))] \bigvee [\Sigma SU(2) \wedge (\mathbb S^6 -SO(3)].$$

\end{exm}

\begin{exm}: $Hom(\mathbb Z^n, SU(2))$:

The first three terms in the decomposition can be determined, and so
$\Sigma Hom(\mathbb Z^n, SU(2))$ splits off a summand of the form
$$\Sigma\bigvee^n G\bigvee
[\Sigma\bigvee^{\binom n 2} (\mathbb S^6 -SO(3))]\bigvee
[\bigvee^{{\binom n 3}}\Sigma SU(2) \wedge (\mathbb S^6
-SO(3))],$$ The remaining stable summands are not yet understood.
\end{exm}

\section{NDR Pairs and Cofibrations}
\label{sec:Cofibrations}

Useful techniques in the proofs of the stable decompositions above
arise from information concerning the geometry of certain Lie
groups, centralizers of elements, and their relationship to the
homotopy theoretic notion of an NDR-pair. Recall Steenrod's
definition \cite{St} that the pair $(X,A)$ is said to be an NDR-pair
provided there exist continuous maps $$u:X \to [0,1],$$ and
$$h:[0,1] \times X \to X$$ which satisfy the following properties.
\begin{enumerate}
  \item $A = u^{-1}(0)$,
  \item $h(0,x) = x$ for all $x$ in $X$,
  \item $h(t,a) = a$ for all $a$ in $A$, and
  \item $h(1,v)$ is in $A$ for all $v$ in $u^{-1}([0,1))$.
\end{enumerate}  Furthermore, if $(X,A)$ is an NDR-pair, then
the natural inclusion $A \to X$ is a cofibration with cofibre given
by the quotient space $X/A$.

Steenrod proves the following Lemma \cite{St} (stated as Lemma 6.3),
and with the extension in the second part due to May \cite{May}
(page 164).

\begin{lem} \label{lem:NDR.pairs}
\begin{enumerate}
  \item Let $(h,u)$ and $(j,v)$ represent $(X,A)$, and $(Y,B)$ as NDR-pairs.
Then $(k,w)$ represents the product pair $$(X,A) \times (Y,B) = (X
\times Y, X \times B \cup A \times Y)$$ as an NDR-pair, where
$w(x,y) = \min(u(x), v(y))$ and
\[
k(t,x,y) =
\begin{cases}
(h(t,x),j([u(x)/v(y)]t,y)) & \text{if $v(y) \geq u(x)$,}\\
(h([v(y)/u(x)]t,x),j(t,y)) & \text{if $u(x) \geq v(y)$.}
\end{cases}
\]

\item Let $(h,u)$ represent $(X,A)$ as an NDR-pair.
Then $(h_j,u_j)$ represents
$$(X,A)^j = ( X^j, \cup_{1 \leq i \leq j}X^{i-1}\times A \times X^{j-i})$$ as an NDR-pair
(which is $\Sigma_j$-equivariant) where $$u_j(x_1,x_2, \cdots, x_j)
= \min(u(x_1),u(x_2), \cdots, u(x_j)),$$ and
$$h_j(t,x_1,x_2, \cdots, x_j) = (h(t_1,x_1), \cdots, h(t_j,x_j))$$
with
\[
t_i =
\begin{cases}
t \cdot \min_{j \neq i}(u(x_j)/u(x_i)) & \text{if some $u(x_j) < u(x_i)$, $j \neq i$,}\\
t & \text{if all $u(x_j) \geq u(x_i)$, $j \neq i$.}
\end{cases}
\]
\end{enumerate}
\end{lem}

This lemma can be used to check that $(Hom(\mathbb Z^n,G),
S_{n}(G))$ is an NDR-pair. However, a modification is required to
prove the full stable splitting in Theorem
\ref{thm:stable.decompositions.for.general.G}. That result requires
that all the $(S_{n}(q,G),S_{n}(q+1,G))$ are NDR-pairs as long as
$S_{n}(q+1,G)$ is non-empty. Thus, let $(X,A)$ denote an NDR-pair
with $F_rX^n$ the subspace of $X^n$ given by $\{(x_1, x_2, \cdots,
x_n)\}$ for which at least $r$ of the $x_j$ are in $A$.

\begin{lem} \label{lem:NDR.pairs.for.fat.wedges}
Let $(h,u)$ represent $(X,A)$ as an NDR-pair. Fix integers $ 1 \leq
r \leq n$. Then $(h_n,u_{r,n})$ represents $(X^n,F_rX^n)$ as a
$\Sigma_n$-equivariant NDR-pair where
$$u_{r,n}(x_1,x_2, \cdots, x_n)
= 1/r \cdot \min_{1 \leq i_1 < i_2< \cdots < i_r \leq n}(
u(x_{i_1})+u(x_{i_2})+ \cdots + u(x_{i_r})),$$ and $$h_n(t,x_1,x_2,
\cdots, x_j) = (h(t_1,x_1), \cdots, h(t_n,x_n))$$ with
\[
t_i =
\begin{cases}
t \cdot \min_{i \neq n}(u(x_m)/u(x_i)) & \text{if some $u(x_m) < u(x_i)$, $m \neq i$,}\\
t & \text{if all $u(x_m) \geq u(x_i)$, $m \neq i$.}
\end{cases}
\]
\end{lem}

\begin{proof}
Notice that
\begin{enumerate}
  \item $u_{r,n}$ is continuous,
  \item $u_{r,n}(x_1,x_2, \cdots, x_n) = 0$ if and only if
$u(x_{i_1})+u(x_{i_2})+ \cdots + u(x_{i_r}) = 0$ for some sequence
$1 \leq i_1 < i_2< \cdots < i_r \leq n$. Since $0 \leq u(x) \leq 1$,
it follows that $u(x_{i_t}) = 0$ for $1 \leq t \leq r $, and so
$x_{i_t}$ is in $A$. Hence $u_{r,n}^{-1}(0) = F_rX^n$.
\end{enumerate}

Finally, notice that

\begin{enumerate}
  \item $h_{n}(-,-)$ is continuous,
  \item $h_n(0,x_1,x_2, \cdots, x_n) = (x_1,x_2, \cdots, x_n)$ as
$h(0,x) = x$,
  \item $h_n(t,x_1,x_2, \cdots, x_n) = (x_1,x_2, \cdots, x_n)$
if $(x_1,x_2, \cdots, x_n)$ is in $F_rX^n$ for all $t$ in $[0,1]$ as
$h(t,x)$ is in $A$ for $x$ in $A$, and
 \item $h_n(1,x_1,x_2, \cdots, x_n)$ is in $F_rX^n$
if $(x_1,x_2, \cdots, x_n)$ is in $u_{r,n}^{-1}([0,1))$.
\end{enumerate}

The lemma follows.
\end{proof}

The next result provides a practical
criterion used in the next section for testing whether certain Lie
groups
have \textsl{cofibrantly commuting elements}.

\begin{thm}\label{thm:NDR.pairs.for.cofibrantly.commuting.G}
Let $(h,u)$ be a representation of $(G, \{1_G\})$ as an NDR-pair
with the additional property that for each $g \neq 1_G$ in $G$ with
$0 \leq t < 1$, the centralizer of $h(t,g)$, $Z(h(t,g))$, equals the
centralizer of $g$, $Z(g)$. If $S_{n}(r,G)$ is non-empty, then
$(h_n,u_{r,n})$ the representation of $(G^n,F_rG^n)$ as a
$\Sigma_n$-equivariant NDR-pair as given in Lemma
\ref{lem:NDR.pairs.for.fat.wedges}, restricts to give a
representation of $(S_{n}(r-1,G), S_{n}(r,G))$ as a
$\Sigma_n$-equivariant NDR pair. Thus $G$
has \textsl{cofibrantly commuting elements}.
\end{thm}

\begin{proof}
The hypothesis of Theorem
\ref{thm:NDR.pairs.for.cofibrantly.commuting.G} is that
$(h,u)$ is a representation of $(G, \{1_G\})$ as an NDR-pair with
the additional property that for each $g \neq 1_G$ in $G$ with $0
\leq t < 1$, the centralizer of $h(t,g)$ equals the centralizer of
$g$.

By Lemma \ref{lem:NDR.pairs.for.fat.wedges}, the pair of maps
$(h_n,u_{r,n})$ represents $(G^n, F_rG^n)$ as an NDR-pair where
$$u_{r,n}(x_1,x_2, \cdots, x_n) = 1/r \cdot \min_{1 \leq i_1 < i_2<
\cdots < i_r \leq n}( u(x_{i_1})+u(x_{i_2})+ \cdots + u(x_{i_r})),$$
and $h_n(t,x_1,x_2, \cdots, x_j) = (h(t_1,x_1), \cdots, h(t_n,x_n))$
with
\[
t_i =
\begin{cases}
t \cdot \min_{i \neq m}(u(x_m)/u(x_i)) & \text{if some $u(x_m) < u(x_i)$, $m \neq i$,}\\
t & \text{if all $u(x_m) \geq u(x_i)$, $m \neq i$.}
\end{cases}
\]

Thus the homotopy defined by $h_n(t,x_1,x_2, \cdots, x_n) =
(h(t_1,x_1), \cdots, h(t_n,x_n))$ satisfies the following
properties:
\begin{enumerate}
\item The centralizer of $h(t,g)$ equals the centralizer of $g$ for all
$g \neq 1_G$ and $0 \leq t < 1$.
\item The element $h_n(t,x_1,x_2, \cdots, x_n)$ is in $Hom(\mathbb
Z^n,G)$.
  \item The homotopy $h_n(t,x_1,x_2, \cdots, x_n)$ preserves the subspaces $S_{n}(r,G)$ and
$S_{n}(r-1,G)$.
\end{enumerate}

Thus the function $h_n$ restricts to
$$h_n: [0,1] \times S_{n}(r,G) \to S_{n}(r,G)$$ and there is a commutative diagram
\[
\begin{CD}
[0,1] \times S_{n}(r,G) @>{h_n}>>  S_{n}(r,G)   \\
@VV{1 \times I_{r}}V             @VV{I_{r}}V          \\
[0,1] \times S_{n}(r-1,G) @>{h_n}>> S_{n}(r-1,G)
\end{CD}
\] for which $I_{r}:S_{n}(r,G) \to S_{n}(r-1,G)$ denotes the natural
inclusion. Since $$u_{r,n}^{-1}(0) \cap S_{n}(r-1,G) =
S_{n}(r,G),$$ the pair $(h_n,u_{r,n})$ restricts to a representation
of $(S_{n}(r-1,G),S_{n}(r,G))$ as an $NDR$ pair. Thus $G$
has cofibrantly commuting elements.
\end{proof}

\section{Closed Subgroups of $GL(n,\mathbb C)$ Have Cofibrantly
Commuting Elements}

In this section it will be shown that any closed subgroup
of $GL(n,\mathbb C)$ has cofibrantly commuting elements,
as mentioned in the introduction (Theorem \ref{thm:closed.subgroups.of.GL.}).
Assume that $G$ is any such group,
with Lie algebra $\mathcal G$. The adjoint representation $$Ad: G
\to Aut(\mathcal G)$$ is defined by the equation $Ad_b(y) = b\cdot
y\cdot b^{-1}$ for $b$ in $G$ and $y$ in $\mathcal G$. Further
assume that $B_{\epsilon}(\vec 0)$ is a ball of radius $\epsilon >
0$ with center at the origin $\vec 0$ in $\mathcal G$ and closure
$\bar B_{\epsilon}(\vec 0)$ such that the restriction of the
exponential map
$$exp:\bar B_{\epsilon}(\vec 0) \to G$$ is a homeomorphism onto its'
image.

\begin{rem} The map $exp$ is a local homeomorphism in case
$G$ is a finite dimensional Lie group \cite{Hum}, but may not be an
embedding on all of $\mathcal G$ in case $G$ is compact.
\end{rem}

\begin{lem} \label{lem:epsilon}
Let $y$ denote an element in $B_{\epsilon}(\vec 0)$ which is not
$\vec 0$. Then the centralizer of $exp(y)$ in $G$ is contained in
the centralizer of $exp(ty)$ in $G$ for any real number $t$.
Furthermore, if $ y \neq \vec 0$ and $ 0 < t \le 1$, the
centralizers of $exp(ty)$ and $exp(y)$ in $G$ are equal.
\end{lem}

\begin{proof}
The following equation is classical \cite{Hum}: $$b\cdot
(exp(y))\cdot b^{-1} = exp(Ad_b(y))$$ where $Ad_b(y) = b\cdot
(y)\cdot b^{-1}$.

Assume that $exp(y)$ and $b$ commute. Thus
$$ exp(y) = b\cdot(exp(y))\cdot b^{-1} = exp(Ad_b(y)).$$
Hence $Ad_b(y)= b\cdot y\cdot b^{-1} = y$ as
$exp$ restricted to $B_{\epsilon}(\vec 0)$ is a homeomorphism. Thus
$b\cdot (ty)\cdot b^{-1} = ty$ for any scalar $t$ as $t$ is central.
Furthermore, $b\cdot(exp(ty))\cdot b^{-1} = exp(Ad_b(ty)) = exp(ty)$
and $exp(ty)$ commutes with $b$. Thus the centralizer of $exp(y)$ is
contained in the centralizer of $exp(ty)$ for all $t$.

Next assume that $ 0 < t \le 1$ and consider the centralizer of
$exp(ty)$. Since $ty$ is in $B_{\epsilon}(\vec 0)$ with $ty \neq
\vec 0$, the analogous argument gives that the centralizer of
$exp(ty)$ is in the centralizer of $exp(sty)$ for all real numbers
$s$. The lemma follows by setting $s=1/t$.

\end{proof}

Define functions $$u,s: G \to [0,1]$$ by the formulas

\[
u(g) =
\begin{cases}
{2|y|/\epsilon} & \text{if $g = exp(y)$ for $g$ in $exp(\bar B_{\epsilon/2}(\vec 0))$,}\\
1  & \text{if $g$ is in $G - exp(B_{\epsilon/2}(\vec 0))$}
\end{cases}
\] and
\[
s(g) =
\begin{cases}
1 & \text{if $g = exp(y)$ for $g$ in $exp(\bar B_{\epsilon/2}(\vec 0))$,}\\
2-u(g) & \text{if $g = exp(y)$ for $g$ in $exp(\bar B_{\epsilon}(\vec 0))- exp(B_{\epsilon/2}(\vec 0))$,}\\
0  & \text{if $g$ is in $G - exp(B_{\epsilon}(\vec 0))$.}
\end{cases}
\] Notice that the functions $s$, and $u$ are both continuous functions.
In addition, $u^{-1}(0) = \{1_G\}$ and $u^{-1}([0,1))=
exp(B_{\epsilon/2}(\vec 0))$.

Next, define a function
$$h: [0,1] \times G \to G$$ by the formula

\[
h(t,g) =
\begin{cases}
exp((1-t)y) & \text{if $g = exp(y)$ for $g$ in $exp(\bar B_{\epsilon\mathfrak /2}(\vec 0))$,}\\
exp((1-s(g)t)\cdot y) & \text{if $g = exp(y)$ for $g$ in $exp(\bar
B_{\epsilon}(\vec 0))- exp(B_{\epsilon\mathfrak /2}(\vec 0))$, and}\\
g  & \text{if $g$ is in $G - B_{\epsilon}(\vec 0)$.}
\end{cases}
\]

\begin{prop} \label{prop:NDR.epsilon}
If $G$ is any closed subgroup of $GL(n, \mathbb C)$, then the pair
$(h,u)$ is a representation of $(G, \{1_G\})$ as an NDR pair.
Furthermore, the representation $(h,u)$ satisfies the additional
condition that for each $g \neq 1_G$ in $G$ with $0 \leq t < 1$, the
centralizer of $h(t,g)$, $Z(h(t,g))$, equals the centralizer of $g$,
$Z(g)$.
\end{prop}

\begin{proof}
Observe that the homotopy $h$ is continuous with $u: G \to [0,1]$.
Furthermore
\begin{enumerate}
  \item $h(0,g) = g$ for all $g$ in $G$,
  \item $h(t,1_G) = 1_G$ for all $t$ in $[0,1]$ and
  \item $u^{-1}([0,1))= exp(B_{\epsilon/2}(\vec 0))$, thus
$h(1,g) = 1_G$ if $u(g)<1$
\end{enumerate} by a direct verification. Thus $(G, \{1_G\})$ is an NDR pair.

Notice that an inspection of the definition gives that the value of
$h(t,g)$ is either
\begin{enumerate}
  \item $g$,
  \item $exp((1-s)y)$ for $y \ne \vec 0$ and $1-s \neq 0$ or
  \item $1_G$.
\end{enumerate} It suffices to consider the centralizer of $h(t,g)$
in the first two cases. Observe that in these cases, the centralizer
of $h(t,g)$ is equal to that of $g$ by Lemma \ref{lem:epsilon}. Thus
the representation $(h,u)$ satisfies the additional condition that
for each $g \neq 1_G$ in $G$ with $0 \leq t < 1$, the centralizer of
$h(t,g)$, $Z(h(t,g))$, equals the centralizer of $g$, $Z(g)$. The
proposition follows.
\end{proof}

\begin{thm}
If $G$ is a closed subgroup of $GL(n, \mathbb C)$, then $G$ has
cofibrantly commuting elements.
\end{thm}

\begin{proof} The previous proposition gives a verification
of the hypotheses of Theorem
\ref{thm:NDR.pairs.for.cofibrantly.commuting.G} thus
implying that if $G$ is a closed subgroup of $GL(n, \mathbb C)$,
then $G$ has cofibrantly commuting elements. The theorem
follows.
\end{proof}

\bibliographystyle{amsalpha}

\begin{thebibliography}{999}

\bibitem{ACC} A.~Adem, D.~Cohen, and F.~R.~Cohen, {\em On representations
and $K$-theory of the braid groups}, Math. Annalen
\textbf{326}(2003) 515--542.

\bibitem {AM} S.~Akbulut and J.~McCarthy, \textbf{Casson's Invariant for
Oriented Homology Spheres}, Mathematical Notes 36, Princeton
University Press (1990).

\bibitem{BFM} A.~Borel, R.~Friedman and J.~Morgan,
{\em Almost commuting elements in compact Lie groups} Mem. Amer.
Math. Soc. 157 (2002), no. 747, x+136 pp.

\bibitem{Cohen-AJM} F.~R.~Cohen, {\em On the mapping class groups for punctured spheres, the
hyperelliptic mapping class groups}, ${\rm SO}(3)$, and ${\rm
Spin}\sp c(3)$. Amer. J. Math. 115 (1993), no. 2, 389--434.

\bibitem{Dold} A.~Dold, \textbf{Lectures on Algebraic Topology},
Springer--Verlag Grundlehren 200 (1980).

\bibitem{FrickeKlein} R.~Fricke and F.~Klein, {\em Vorlesungen \"uber die Theorie
der automorphen Funktionen}, Bd. 1, Teubner, Leipzig, 1897, 2nd
edition, idem, 1926.

\bibitem{Goldman} W.M.~Goldman, {\em Topological components of the space of
representations}, Invent. Math. \textbf{93}(1988), no. 3, 557-607.

\bibitem{Hum} J.~Humphries, \textbf{Introduction to Lie Algebras and Representation
Theory}, Second printing, revised. Graduate Texts in Mathematics, 9.
Springer-Verlag, New York, 1978.

\bibitem{HL} N.~Huo and C.C.~Liu, {\em Connected Components of
Surface Group Representations}, Int. Math. Res. Not. \textbf{44}
(2003), 2359--2372.

\bibitem{Kac-Smilga} V.G.~Kac and A.V.~Smilga,
{\em Vacuum structure in supersymmetric Yang-Mills theories with any
gauge group},  The many faces of the superworld,  185--234, World
Sci. Publishing, River Edge, NJ, 2000.


\bibitem{KapMill} M.~Kapovich, and J.J.~Millson, {\em On representation varieties
of Artin groups, projective arrangements and the fundamental groups
of smooth complex algebraic varieties}, Publications Mathematiques
IHES, 88 (1998), pp. 5--95.

\bibitem{Lannes} J.~Lannes, {\em Th\'eorie homotopique des groupes de
Lie (d'apr\'es W.~G.~Dwyer et C.~W.~Wilkerson)}, S\'eminaire
Bourbaki, Vol. 1993/1994. Ast\'erisque (227 1995), Exp. No. 776, 3,
21-45.

\bibitem{Li} J.~Li, {\em The space of surface group representations},
Manuscript. Math. \textbf{78}(1993), no. 3, 223-243.

\bibitem{May} J.P.~May, {\em The geometry of iterated loop spaces}
Lecture Notes in Mathematics \textbf{171}, 1972.

\bibitem{St} N.~E.~Steenrod, {\em A convenient category of
topological spaces}, Michigan Math. Journal \textbf{14} (1967),
133-152.

\bibitem{W} B.~Wajnryb,
{\em A simple presentation for the mapping class group of an
orientable surface}, Israel J. Math. 45 (1983), no. 2-3, 157--174.

\end{thebibliography}

\end{document}